\documentclass[12pt]{article}

\usepackage{graphicx}
\usepackage{url}

\newcommand{\beq}{\begin{equation}}
\newcommand{\eeq}{\end{equation}}

\title{ Euphemia Lofton Haynes:
 \\ Bringing Education Closer to the ``Goal of Perfection''  \footnote{This phrase is from Haynes' {\em Mathematics-Symbolic Logic} 1945 address \cite{HaynesLogicQuote}. }}
   \author{Susan E. Kelly  \and Carly Shinners \and Katherine Zoroufy \thanks{Susan Kelly (skelly@uwlax.edu) is a professor of mathematics at the University of Wisconsin La Crosse.  Carly Shinners is currently teaching high school mathematics in South Korea and Katherine Zoroufy works at Epic in Madison Wisconsin.  Shinners and Zoroufy began their work on this project as an undergraduate research project directed by Kelly.} }
     \date{ }

\begin{document}

\maketitle

\begin{figure}
\centering
\includegraphics[height=100mm]{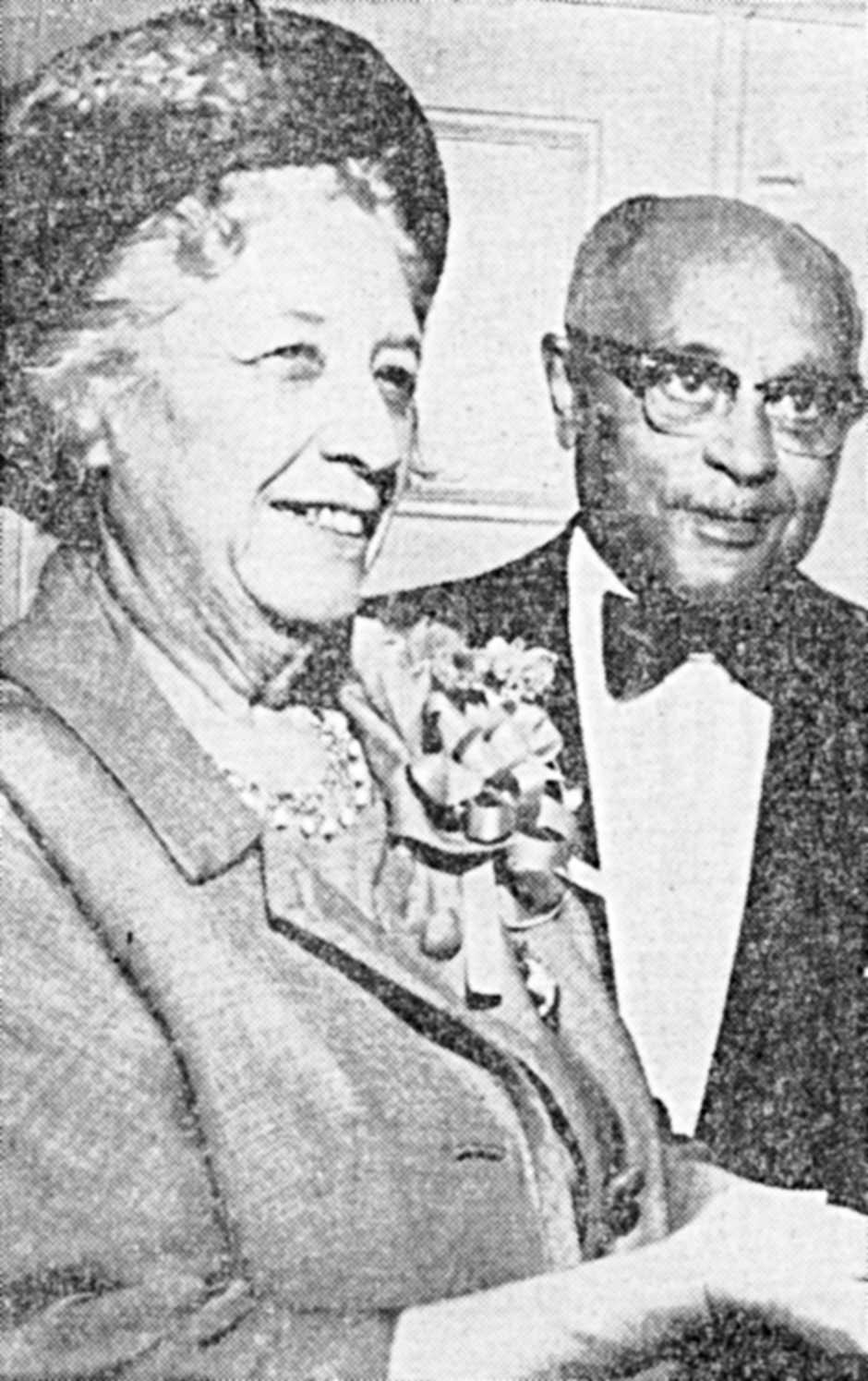}\\
\caption{Haynes was named Lady of the Year by the DC ``Oldest Inhabitants'' in 1967.}
\end{figure}

 Martha Euphemia Lofton Haynes was the first African American woman to receive a PhD in mathematics.  She grew up in Washington DC, earned a bachelors degree in mathematics from Smith College in 1914, a masters in education from University of Chicago in 1930, and a doctorate in mathematics from the Catholic University of America in 1943. Haynes spent over forty-five years teaching in Washington DC from elementary and secondary level to university level. She was active in many community service organizations where she served in leadership roles and received numerous honors including being named a fellow of the American Association for the Advancement of Science and being awarded a Papal Medal. She was a member of the Washington DC school board from 1960 - 1968, serving as president from June 1966 through July 1967.  She played a leadership role in ending the tracking system, which she argued discriminated against African American students by assigning them to education tracks that did not prepare them for college.  This fight culminated in the 1967 Hobson v Hansen court case, in which the judge ruled that tracking was discriminatory towards poor and minority students.


\section{Setting the Stage}

Seldom does the story of one's life begin at birth. In documenting the life and accomplishments of Haynes, it is important to also look at some of the history which sets the stage for her life.  Haynes was born on September 11, 1890 in Washington DC \cite{EuphemiaDeathRecord}, just twenty-five years after the end of the United States Civil War. This location, the time period, and her race played major roles in shaping her life.

Washington DC by its very nature has always been unique. This is illustrated in its history related to race.  The Declaration of Independence states ``We hold these truths to be self-evident, that all men are created equal... .'' Yet, slavery was allowed to continue. George Washington was a slave owner. The White House and Capital Building were built with slaves.  Philip Reid, a slave, worked on the {\em Statue of Freedom}, which sits atop the Capital dome \cite{WHCap, Reid}.  Slaves, some likely owned by Martha Washington, quarried the sandstone used to build the Smithsonian Castle on the National Mall \cite{Castle}. In 1800, one fourth of DC's population was African American, most of whom were slaves.  By 1830, the majority of African Americans in the capital were free. While free African Americans had access to education in Washington DC, in neighboring states such as Virginia such education was against the law.  Also, at this time, public auctions of slaves continued in the nation's capital.  In April 1862, nine months prior to the general Emancipation Proclamation, President Lincoln signed the DC Compensation Emancipation Act, which began the process of freeing slaves in DC while compensating their owners for their deemed value \cite{McQuirter, DChistory}.

\begin{figure}

\centering
\includegraphics[height=150mm]{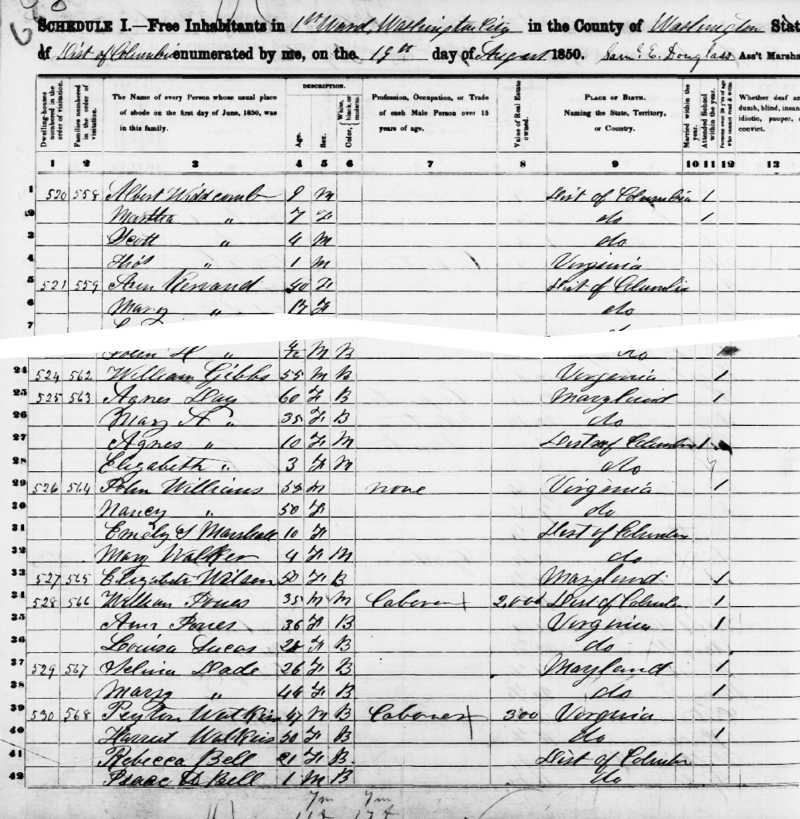}\\
\caption{Euphemia's Great-Great Grandmother Agnes Day, line 25, Great-Grandmother Mary and Grandmother Agnes are listed as Free Inhabitants of Washington D.~C. in the 1850 Census \cite{census1860Day}.}
\end{figure}

Haynes' mother, (Anne) Lavinia Day Lofton, was a native of Washington DC and taught kindergarten in the public schools. She was active in her Catholic church as an organist, children's choir director, and Sunday school teacher \cite{HLBiographicalNotes}. Haynes' mother's side of the family can be traced back to Haynes' great-great grandmother Agnes Day who was born in Maryland in about 1790. In the 1850s and 1860s, the Day family members were listed as ``Free Inhabitants'' of Washington DC.  Census data from 1850-1870 list no adult males living in the household.  Adult females were listed with occupations of domestic servant, housekeeper, and washing \cite{census1850Day,census1860Day,census1870Day}. Little information was found on Haynes' maternal grandfather, Charles Delany, sometimes spelled Delaney, who was born in Maryland. Charles' mother, Grace Delany, was born in Maryland, and in the 1870 and 1880 censuses, she was living in Washington DC and is listed as a servant/housekeeper of Thomas Abbot, a white man \cite{LaveniaDeath,census1870Delany,census1880Delany}. Both Grace Delany and Haynes' maternal grandmother, Agnes Day, were listed as mixed race \cite{census1870Delany,census1860Day}.

Haynes' father, William S. Lofton, was born in the 1860's in Batesville, Arkansas and moved with his parents to Washington DC prior  to the 1870 census\cite{census1870Lofton}. William was a graduate of Howard University and became a successful dentist and a member of the Board of Directors of the Capital Savings Bank. He was a Catholic lay leader on the national level, who pressed the Church to establish schools for African American children and worked to help create Catholic organizations for African Americans in the 1890s, when racial prejudices were increasing \cite{Aristocrates, HLBiographicalNotes, Moore}.   In the 1870's census, William's father, also named William, was listed as a laborer born in Kentucky who could not read or write, and his mother, Martha, was a housekeeper who was born in Maryland or Missouri who could not write \cite{census1870Lofton,census1900Lofton}. By the 1900 census, Martha was listed as a housekeeper with the skills to both read and write. Haynes' father was listed as being mixed race according to the 1910 census \cite{census1900Lofton, census1910Lofton}.

Although no genealogical records prior to 1870 were found for the Lofton family, one can consider some possibilities of their lives based on the history of the time. One major obstacle in searching African American families prior to 1870 is the fact that slaves were not counted in census data as named individuals.  Often, slave records only record numbers of slaves and details such as age and gender.  It is not know if the Lofton family was free or held in slavery, so we can only consider likelihoods.  Arkansas, the state of Haynes' father's birth, became a state in 1836.  During it's territory years and early years of statehood, some free blacks lived in communities, often in more remote areas of Arkansas.  For instance, north-west of Batesville, on the border with Missouri, there was a significant sized free black community.  However, in 1859, the law {\em An Act to Remove the Free Negroes and Mulattoes from the State} was passed which threatened slavery to such individuals who resided in Arkansas after January 1, 1860.  Although the law was repealed in 1863, in 1860, Arkansas was the slave state with the fewest free blacks \cite{Higgins}.

By the time Haynes was born, slavery had been abolished; however, opportunities and rights of black Americans were not equal to those of whites. In the United States in 1890, at the time of her birth, 45\% of blacks fourteen years old or above were illiterate, while only 6\% of whites of that age were.\cite{Literacy} Haynes' parents could read and write, but as stated earlier, some of her grandparents could not \cite{census1900Lofton,census1900Euphemia}. Schools in DC were segregated along with many other aspects of society. Haynes lived through the Civil Rights years including the 1954 Brown v Board of Education decision, which ruled that segregation of public schools was illegal.  Her background and location set the stage for her fight for equal education opportunities in the nation's capital.


\section{Early Life and Family}

Martha Euphemia Lofton Haynes prefered being called Euphemia rather than Martha. Her father, William Lofton, was a member of the Washington ``black 400,'' a small group of fewer than 100 families in Washington DC who were considered aristocrats of color, a distinction often based on family background, occupation, color and generations removed from slavery \cite{Aristocrates}. Because Haynes' father was a member of this elite group and her mother was a school teacher, Haynes had opportunities that most African Americans in DC would not have at this time.  However, her family situation began to experience turmoil early in her life.

Haynes' family began to break up during her early youth.  Her brother Joseph was born in 1893 \cite{HLBiographicalNotes}, but by 1895 her parents had separated.  In a letter William sent to Lavinia in December of 1895, William stated that he had not seen their children except by accident since September of that year and accused his wife of intentionally keeping the children from him.  He told her that he had prepared a comfortable home for her and their children and also inquired to let him know if the children were in need of shoes \cite{WilliamLetterDec1895}.  In 1900, the divorce was finalized, and custody was given to Haynes' mother, while Haynes' father was given weekly visitation rights \cite{divorce}.  The relationship between Haynes and her father appeared to be strained throughout the rest of William's life, based on letters William wrote to his daughter \cite{William1903, William1904,William1911}.  William Lofton died in 1919. His letters to Haynes were saved, and William Lofton included both of his children in his final will \cite{CertificateTitle, HLBiographicalNotes}.

\begin{figure}
\centering
\includegraphics[height=80mm]{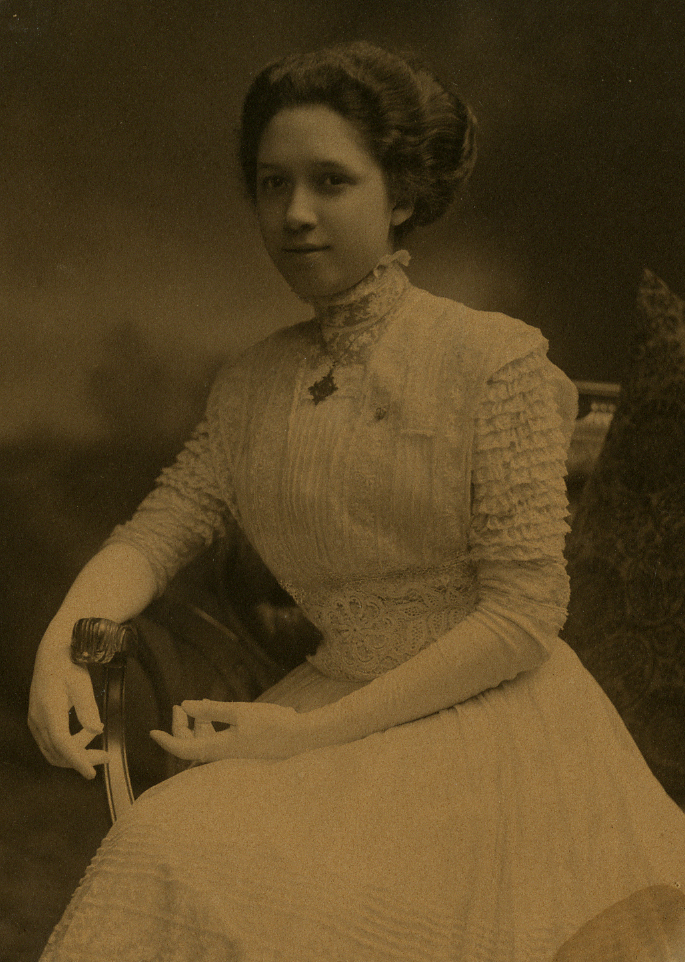}\\
\caption{Euphemia Lofton before marrying Harold Haynes in 1917}
\end{figure}

By the time of the 1900 Census, nine year-old Haynes was living with her family in the home of her uncle and aunt, Benjamin and Anna Swann.  Benjamin was a butler and Anna, Lavinia's sister, was a dressmaker. Other adults in the home were Haynes' mother Lavinia Lofton, who worked as a school teacher, Haynes' grandmother Agnes Day, and Haynes' great-grandmother Grace Delany. Children in the home were Haynes, her seven-year-old brother Joseph, and her two year-old cousin Daniel Swann \cite{census1900Euphemia}. Archive records indicate family struggles at times, but also showed a family that stayed close.  Later in life, Haynes traveled to Chicago frequently when her brother was dying, and in a 1972 interview Haynes spoke of her mother, saying, ``My mother was so successful as a mother because she believed in me'' \cite{change}.

\begin{figure}
\centering
\includegraphics[height=80mm]{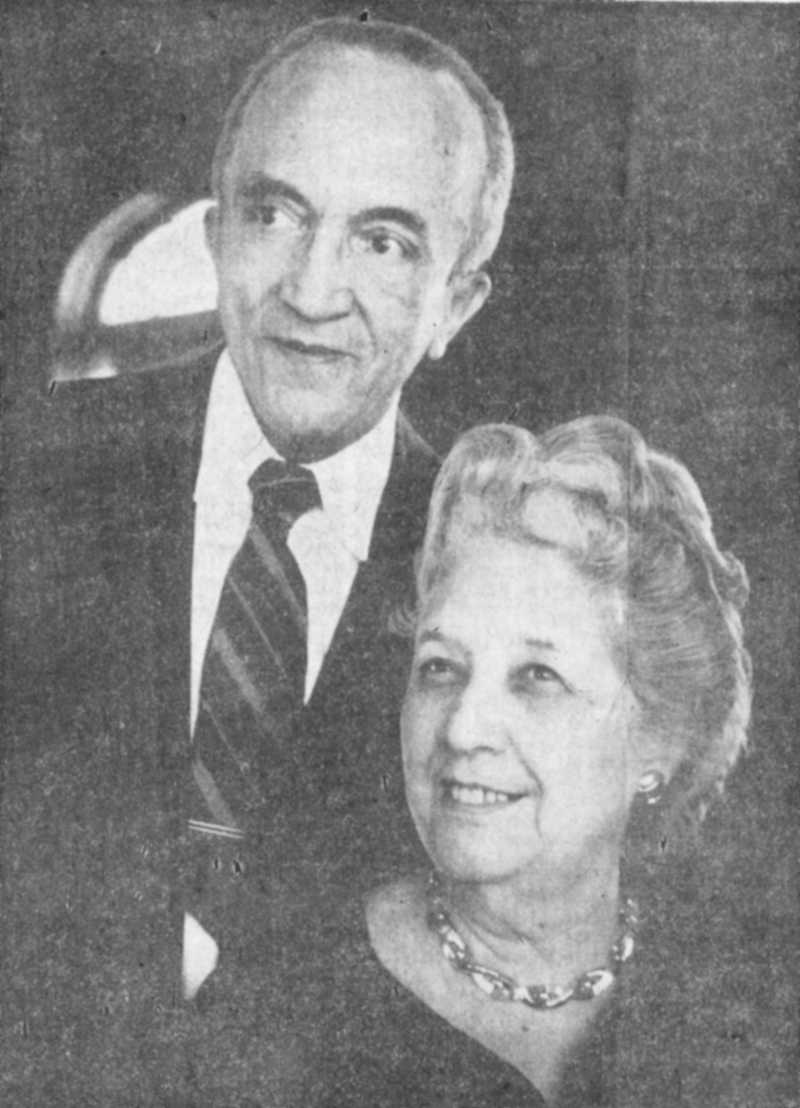}\\
\caption{Euphemia and Harold Haynes}
\end{figure}

Euphemia Lofton married Harold Haynes in 1917. They had known each other as teenagers growing up in the same neighborhood. In 1908, Harold wrote to Euphemia telling her not to worry about her mom wishing to separate them.  He said it made him ``more determined to stick'' and that he had learned ``that it pays to wait for some things.'' Harold graduated from M Street High School in 1906, a year earlier than Euphemia. He earned an electrical engineering degree from the University of Pennsylvania in 1910, a masters in education from the University of Chicago in 1930, and a doctorate in education from New York University in 1946. His teaching career included time at Howard University before moving to teach in DC public schools.  He followed this by serving as a public school principal and then an associate superintendent. He became the superintendent of Black Schools in DC's segregated schools in 1951 and moved to the role of deputy superintendent when schools were desegregated.  He retired from this role in 1958 \cite{HLBiographicalNotes, HaroldLetters,HaroldQual}.

\section {Education}

Haynes graduated as valedictorian from the M Street High School in 1907 \cite{HLBiographicalNotes}. The school's origins can be traced to 1870 when a school for blacks was established in a church basement after the United States Congress defeated a bill for an integrated public school system in the nation's capital. A compromise was a promise of equal quality segregated education. A school building was first constructed in 1890-1891 and was one of the nation's first high schools for African American students. Because there were few opportunities for black professionals, the school was able to obtain teachers whose educational backgrounds surpassed those of most teachers at high schools for white students.  Many of the graduates of M Street High School went on to college or university and became black leaders in DC and beyond \cite{MStreet}. In Haynes' valedictorian speech, her words reflect her actions throughout her life: ``For a person of intelligence is well equipped to solve the problems of life...We must have some defined aim in life and be able to fill competently that position in which we may find ourselves...Let each defeat be a source of a new endeavor and each victory the strengthening of our spirit of gratitude and charity towards the unsuccessful'' \cite{ValSpeech}. As will be addressed in Section 6, lack of funding directed to black schools in DC allowed the schools' quality to greatly diminish, and Haynes would become a voice challenging the unequal education offered in public schools based on race.

Haynes next attended Miner Normal School and graduated in 1909 \cite{HLBiographicalNotes}. This school also has a fascinating history. Myrtilla Miner was born on March 4, 1815 in Brookfield New York to a family with little education.  With her strong desire for a good education, she borrowed and read every book she could find.  With little financial means, she asked to be admitted to school in Clinton, New York and delay payment for a year so that she could find a teaching job to earn the money for the tuition. Upon graduating, she first taught in New York and Rhode Island and later in Mississippi.  In Mississippi, she taught at a school for plantation daughters.  She was horrified with the system of slavery.  She first hoped to assist in plans to free some slaves, but that didn't materialize.  She next asked to teach slave children, but was informed that teaching slaves to read was illegal. She returned home two years later due to health issues and promised herself that when her health allowed, she would set out to educate black children.  With a \$100 contribution and school supplies obtained by begging her friends, she started a school for young black girls on December 3, 1851 in Washington DC.  The school started with six students and grew to forty by the second month. The school faced hostilities from the community including verbal abuse of students, threats of violence, stones thrown through school windows, and the school being set on fire one evening. To help defend the school, Myrtilla practiced shooting a pistol \cite{Miner}. In 1955, after the Brown v Board of Education decision, Miner's school merged with a school for white girls to become the District of Columbia Teacher's College \cite{UDChistory}.

She began teaching elementary school and enrolled in Smith College.  She earned a bachelor of arts with a major in mathematics and a minor in psychology in 1914 \cite{RetirmentNotes}. In letters written to her mother, vivid images of this time period in her life are exposed.  In one such letter, she writes about her determination and perseverance: ``In haste to give you the news that I passed the exam ... I suffered everything from 2-5:30. When I came out I could hardly walk home... When I sat down to take it, he did not want to give it to me and said he knew I would not pass it and it would be fatal for me if I did not. Now that is what I had to brave. I knew if I did not give an excellent paper he would never let me through''\cite{SmithLetters}.

In that same letter, she describes a more pleasant experience: ``Did I tell you too that last week we had to write an original composition in music. Mine was the best in the class and Mr Moog played it before the class.'' In another letter, she wrote a rather humorous account: ``I weigh 122. All my clothes are too tight. When I put on a light dress at dinner now I have to take off flannel petticoat. Of course I don't go out that way.  Please tell Miss Helen and ask her to send me something. Tell her I simply can't get along. One of my petticoats split wide open on me. When I sit down the girls say they are in mortal dread that my dress will split. I had no idea I would get so fat'' \cite{SmithLetters}.

Finally, her letters also express her excitement of learning and of soaking up ideas. In one example, she describes an occasion when she heard Dr. Lyman Abbott, a nationally known Congregationalist pastor and editor of several magazines and periodicals, speak on campus: ``He spoke on the commerce of ideas, the commerce of thought... He just filled us with it. He started  in a way that seemed to shock everyone by saying that we were not going out into the world to serve as it was so often put. Then he went on to explain that it was our business to take as well as give. He said `There is no living being in the world who does not know more about something than anyone of you. From every person even a beggar we can learn something.' Then he said we must remember that what we got in college it was our duty to give to our friends who did not have that opportunity...''\cite{SmithLetters}.

Haynes' graduate work took her to the University of Chicago.  The university was a leader in admitting women and students of color for advanced degrees, with 45 African Americans earning PhDs between 1870 and 1940, the highest number for any institution in the country. Reasons cited for this leadership role include the university being more open than more established institutions in experimenting to establish its identity, its proximity to an established black community, and the ``courage and conviction of a few faculty and administrators, combined with a handful of African American students determined to pursue their education'' \cite{Chicago}.

Haynes earned a masters degree in education from Chicago in 1930 with her thesis {\em The Historical Development of Tests in Elementary and Secondary Mathematics}.  In this work, she traced the development of standardized testing tools from 1900 to 1930. In her conclusions, she noted the difficulty of early tests in determining both the nature and causes of variations in scores. She also discussed later tendencies to use tests to measure progress in learning rather than in classifying pupils \cite{thesism}. These observations would serve her well later in her career as she challenged the tracking system of DC public schools. Euphemia then began taking graduate level mathematics classes from the University of Chicago \cite{HLBiographicalNotes,Vitea, CUAtranscript}.

\begin{figure}
\centering
\includegraphics[height=80mm]{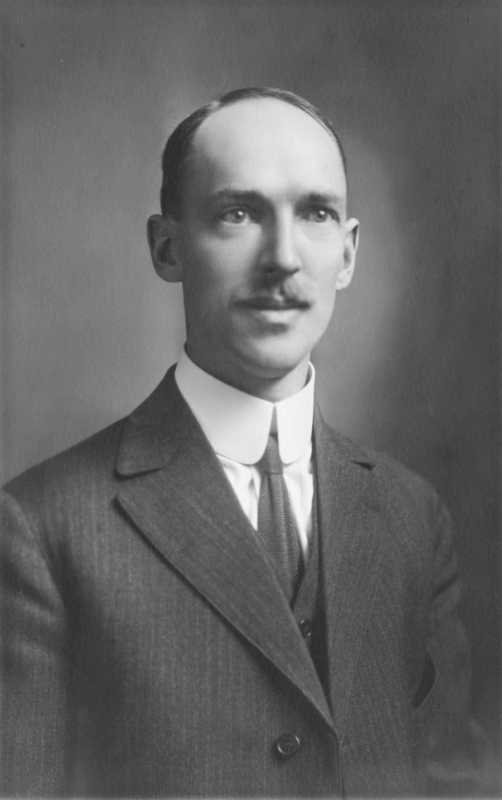}\\
\caption{Haynes' doctoral thesis advisor Aubrey E. Landry was a leading advisor for women.}
\end{figure}

Later, Haynes returned to school to pursue a doctoral degree in mathematics from The Catholic University of America in Washington DC. The university, which attracted many nuns, had produced the fourth highest number of female doctoral students prior to 1940 in the country. Euphemia's thesis advisor, Aubrey Landy, was responsible for directing the mathematical dissertations for all the nuns at the time. He had earned a PhD from Johns Hopkins in 1907 under the direction of Frank Morley, who also was a leading advisor for women in mathematics. Landy's research area was algebraic geometry, and he typically chose thesis topics outside mainstream research for students looking to gain credentials for teaching rather than research\cite{HLBiographicalNotes,Green}.

Haynes' dissertation, {\em Determination of Sets of Independent Conditions Characterizing Certain Special Cases of Symmetric Correspondences}, expanded on work done by another of Landy's students whose thesis was finished in 1938. In Euphemia's work, she examined two ways geometric representations were  defined on parametric rational plane curves in the 1938 thesis and investigated their differences \cite{thesisdr,thesisGalvin}. Like many of Landry's students, Euphemia did not continue to work in this area after graduating, and much of the vocabulary of the thesis is not commonly used today.

\begin{figure}
\centering
\includegraphics[height=120mm]{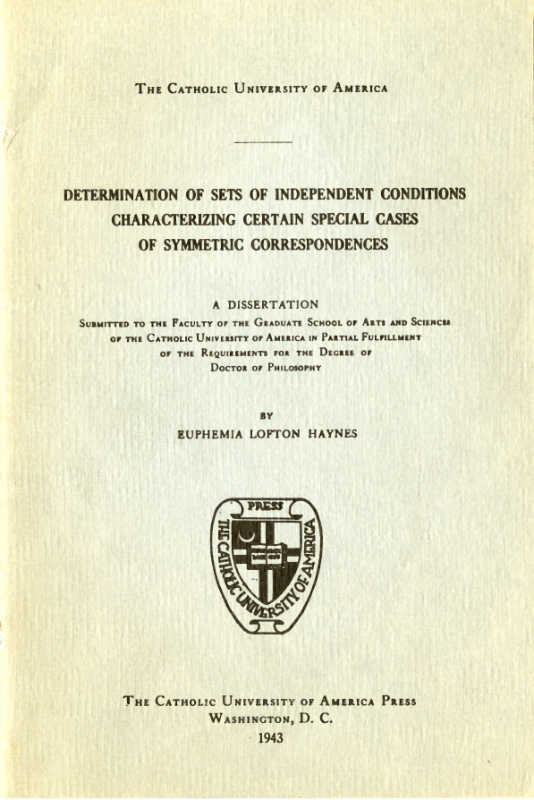}  \\
\caption{Haynes' doctoral thesis studied rational plane curves.}
\end{figure}

Haynes completed her doctoral degree in 1943. She applied her educational experiences to both her teaching career and to the many service activities that she pursued throughout her life.

\begin{figure}
\centering
\includegraphics[height=80mm]{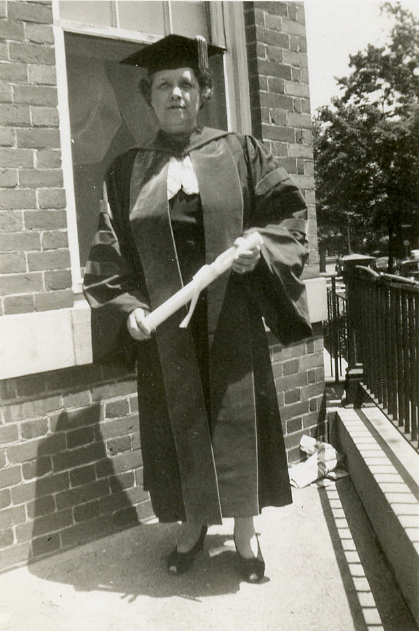}\\
\caption{Haynes received her doctoral degree from The Catholic University of America in 1943.}
\end{figure}


\section{Teaching}

In a 1945 address, {\em Mathematics - Symbolic Logic}, given to junior high and high school mathematics teachers, Haynes eloquently described the full beauty of mathematics framed with logic and the need for teachers to convey this understanding. She stressed the need to devote significant time for observation and reflection to establish truth rather than repetition to cement in facts.  She stated that if not taught correctly, the true nature of mathematics would not be seen: ``Mathematics is no more the art of reckoning and computation than architecture is the art of making bricks, no more than painting is the art of mixing colors.''  She went further in saying, ``...what is the mathematician doing? He is building notions or ideas, he is constructing, inventing, adding to his body of science.  With what is he working? Ideas, relationships, implications, etc. What are his methods? Observations, experimentation, incomplete induction. He is deliberately providing time for reflection and contemplation.''  To be able to teach well, she emphasised the need for teachers to have done serious study into the breadth of the discipline \cite{HaynesLogicQuote}.

Haynes began her teaching career as an elementary teacher in 1909 after graduating from Miner Normal School. She later taught mathematics at various DC public high schools, including becoming the head of the mathematics department at Dunbar High School. She served as professor of mathematics at Miner Teaching College, organizing the department in 1930 and also chairing the department. With the merger of schools after desegregation, she continued as professor and chair of the DC Teacher College until her retirement in 1959 \cite{Vitea}. During part of that time period, she also was a part-time instructor at Howard University \cite{Griffith}.

In a 1960 address on International Communication, Haynes described the quality teacher: ``We remember that professor who was really able to communicate his point of view. Why will we remember him? His vitality, his zeal for truth, the apparent joy and satisfaction he experienced from his endeavors - all these evidences of enthusiasm stimulated us. His enthusiasm was contagious. This great teacher was a great person operating thru the medium of subject matter'' \cite{Communication}.

It seems fitting to close this section on teaching with words she gave at a 1960 high school commencement address: ``I believe there are two requisites for success in life: (1) that one be always a student and (2) that he dedicate himself to the service of others'' \cite{RooseveltHigh}.

\begin{figure}
\centering
\includegraphics[height=60mm]{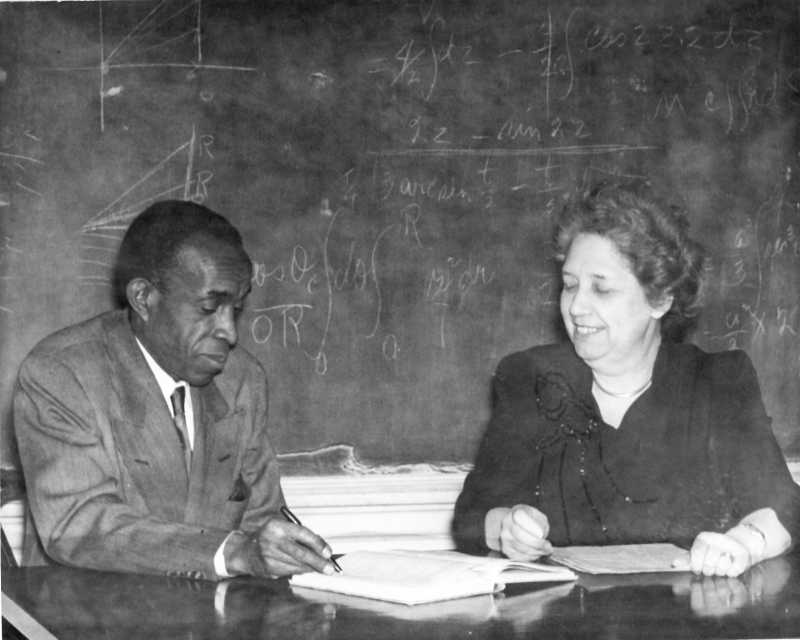}\\
\caption{Euphemia Lofton Haynes in Undated Photo}
\end{figure}


\section{Service and Honors}

Haynes lived according to her message for service. A complete list of her activities is likely not possible, but a sample gives a picture of her tireless dedication towards causes she valued.

She started service early in her life by following in her mother's footsteps as organist and choir member in her church. She continued church-related work with service on the board of Catholic Charities. She was also on the executive committee of the  National Conference of Christians and Jews and was one of the founders of the Catholic Interracial Council. As a member of the the Washington Archdiocesan Council of Catholic Women, her service included time as a vice-president and as president. On the national level, she served as president for fourteen years of the New Federated Colored Catholics of America.  In 1959, she was awarded the Pro Ecclesia et Pontifice, the highest medal that the Pope can award to laity \cite{Vitea, Pope}.

In an address she delivered in the 1960s entitled ``What Faith Means to Me'', Haynes stated, ``My faith also teaches me that every human being is created in the image and likeness of God. .. As one moves about his daily work he influences the lives of his brothers. It is his obligation to be certain therefore that his influence contributes always to the salvation of these souls. An outstanding example of this responsibility in action is the `living wage' for workers ... it is a concrete and material expression of man's consciousness of his responsibility to his fellow man in the field of labor relations.'' She went further to speak of the National Labor Relations Act, which guarantees employees in the private sector the right to collective bargaining, and the Social Security Act as additional examples of policies which modeled these principles in practice \cite{Faith}.

Her application of her faith to serving others and to social justice issues is reflected in many areas of service. She was a member of the United Service Organization (USO) from  1943 to 1965, including twelve years on the national board of directors and service on the National Committee on Service to Negroes and on the National Committee on Service to Women and Girls \cite{Vitea}. In a letter from 1946, she wrote ``Like many of the volunteer workers with the USO I am deeply appreciative of the opportunity afforded me because I realize that I have gained far more than I have given'' \cite{USOletter}.  She also served on the National Committee of the Girl Scouts and on the local level was chairman of the advisory board of the Fides Neighborhood House, which provided activities for children and distributed food, clothing, and other basic needs to families \cite{Vitea, Fides}. She also accepted an invitation to serve as a delegate to the 1960 White House Conference on Children and Youth \cite{WH}.

\begin{figure}
\centering
\includegraphics[height=90mm]{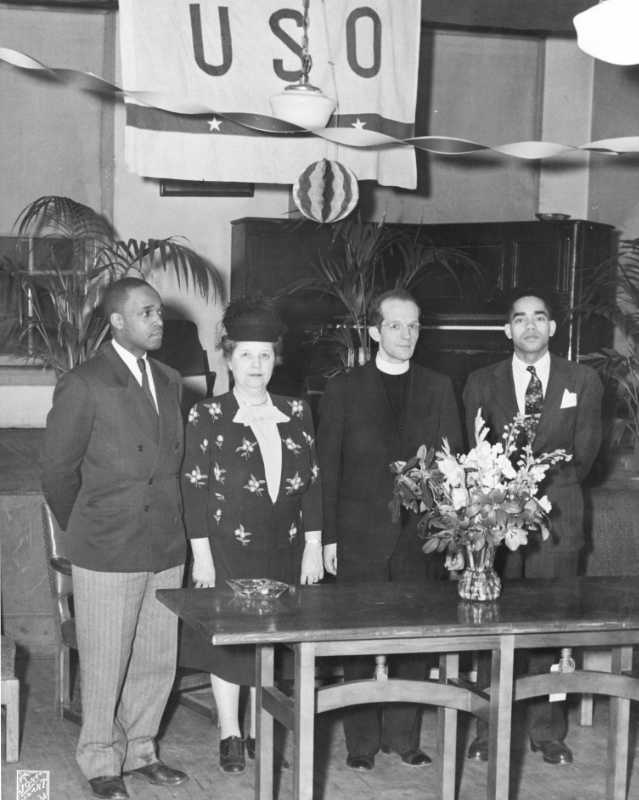}\\
\caption{Haynes was a member of the USO from 1943 to 1965.}
\end{figure}

Haynes was a member of the American Mathematical Society, a Fellow of the American Association for the Advancement of Science, president of the National Association of College Women - DC Branch, and chairman of the Committee on Education for Sigma Delta Epsilon, a graduate fraternity of women scientists \cite{AMS, Vitea,SDE}.

Some of the distinguished awards related to her residence in the nation's capital included an invitation to a reception from the Ambassador of Ethiopia \cite{Etho}, a dinner invitation in honor of Mrs. Eleanor Roosevelt as United States Delegate to the General Assembly of the United Nations \cite{Roosevelt}, and being given the ``Lady of the Year'' award in 1967 by the Oldest Inhabitants Inc., an organization that originally formed in 1912 for black residence who had lived in the District for a set number of years.  Vice-President Herbert Humphrey sent a letter to Haynes to congratulate her on this honor \cite{LadyYr, Aristocrates,VP} .

\section{D.C. School Board and Battle for Improved Integration}

In 1960, shortly after retiring from a 47-year teaching career, Haynes was an invited panelist discussing how retirement could be viewed as a new career.  Haynes' personal notes for her discussion state, ``Concept of retirement  - Not as a termination, but as a new challenge or opportunity. Confidence of Self: In life one has accumulated valuable experiences and new talents.  Retirement offers opportunities to use these in new and different areas'' \cite{retirementquote}.  After only eight months of retirement, Haynes was selected to serve on the DC school board to finish the remaining three month term of a member who resigned due to ill health \cite{startSchoolBoard}. Thus, Haynes began a new career, and used it to achieve possibly her greatest impact on education.

To understand Haynes' role in DC public education, it is important to understand some of the history of the DC public school system. As stated earlier, when M Street school was established in the 1800s as the first black high-school in the nation, the quality of education at M Street equalled or surpassed the quality at white schools. However, as the black population of DC grew, new schools were not built, and overcrowding became an issue in black schools.

Carl Hansen left the desegregated school system of Omaha, Nebraska to become the assistant superintendent of the segregated public school system of DC in 1947 \cite{HansenObit}. He wrote that the ``colored branch of the school system was only a dimly outlined appendage of the white structure,'' and he wanted to do all he could to desegregate the schools.  Yet, as part of his job, he enforced segregation while it was the law. For instance, he wrote of the first school sit-in in 1952 when a black minister sat quietly with his young son in a white kindergarten classroom. The man wished to enroll his son in the school that was only blocks from his home.  There was no law against trespassing, and thus he could only be removed if he was disturbing the peace.  The next day, Hansen arrived at the school prior to the man arriving again with his son.  Hansen blocked the entry so that the only way the man could enter would be to forcibly push Hansen.  Police were at the school in such a case to arrest him. Hansen wrote that the man hesitated and then sadly left with his son in hand. Hansen later wrote, ``It is not necessary to hark back to days of slavery to talk about harm to Negroes.  Today almost as much as in days past, Negroes still suffer from the cruelty of attitudes and actions generated by the color of their skins rather than by the character of their behavior'' \cite{Hansen}.

In 1952, twelve year-old Spottswood Bolling and four other students were plaintiffs in the DC Bolling v Sharpe case that argued that black students in DC were attending overcrowded schools and denied access to better white schools.  The case was decided in 1954 by a Federal district judge who ruled against the plaintiffs, but the US Supreme Court decided to review it, combining it with other segregation cases. In 1954, after the Supreme Court's Brown v Board of Education decision, DC schools were desegregated.  Hansen called it ``a miracle of social adjustment unprecedented''\cite{HansenObit}.

Yet desegregation did not solve all problems. Students were required to attend their neighborhood school except for the greatest of reasons. This created largely homogeneous schools based on race and income. Exceptions to the neighborhood school policy tended to be made for affluent white students wishing to leave a predominantly black school \cite{Hansen, HvH}. In addition, after desegregation ended, a four-track curriculum was introduced in 1959 for students at all levels. Based on IQ test scores or opinion of a principal or teacher, students were tracked into {\em Honors} or {\em Regular College Preparatory}, which prepared students for college; {\em General Curriculum}, which educated students for blue-collar type work; or {\em Basic Curriculum} for students deemed to be academically delayed \cite{Tracks, SB1965,1963Statement}.  Hansen had become superintendent in 1958 and was a strong supporter of both neighborhood schools and the track system \cite{HansenObit}.  It was in this setting that Haynes joined the School Board.

In November 1963, Haynes spoke of the lack of validity of IQ tests and the question of whether they measured cause or effect, an area related to her masters research. She pointed out the segregation created by the tracking system and the increase in the number of student dropouts. She stated how ``man made predictions determine the future of the child'' and removed ``freedom of choice.'' She stated that the tracking system was in ``direct opposition to the American ideal ... a free society capable of self-direction and appreciative of the dignity and potential worth of all members.'' She pointed out how it gave those placed in the highest tracks a ``better than thou attitude'' and gave those placed in the lower tracks a ``feeling of inadequacy'' \cite{1963Statement}.

In the December 1963 meeting of the school board, Haynes stated, ``Where in our track system do we find any self-participation in the determination of the kind of educational experience one may have? The importance of the opportunity to succeed has been recognized.  But equally important is the opportunity to fail.  Only here does one meet the challenge to modify his chosen path. Never on the basis of dictation but only thru his right to try and his right to fail can any student accept an evaluation of self so essential to an attitude of self respect, so necessary to the very fundamental appreciation that for every human being there are areas in which he will succeed and there are areas in which he will fail.'' She requested that the superintendent and staff study the issue and replace the track system with a system that meets the needs of all students \cite{EHDec1963}.

In 1964, Haynes spoke again before the school board claiming that the tracking system ``attempts in apartheid-like fashion to separate the underprivileged.'' She said there was little opportunity for students to move between tracks; some tracks did not prepare students for further education, and diplomas from the various tracks carried different values.  She added that, ``a school experience which insures no contact of `my' group with `that' group and preserves the attitude of `we' and `they' cannot lead to a unified citizenry, working towards the same goals.'' The board submitted Haynes' recommendations concerning these issues to Superintendent Hansen \cite{1964Paper, EH1964}.

In early 1965, letters from parents, administrators, and education experts were collected both supporting and denouncing tracking \cite{Letters1965}. In April, a special meeting of the school board was held to consider tracking, and Haynes moved that the track system be abolished
 \cite{SB1965}. In August, Superintendent Hansen agreed to have an outside agency evaluate the system even though Haynes objected saying that money should not be spent to obtain findings that were already known \cite{Afro1965}. Superintendent Hansen hired Harry Passow, a Professor of Education at Columbia University to conduct a study on the situation of the DC public schools.  In October, findings were reported.  It was determined that almost all the students in the {\em Honor} track were white, while almost all of the students in the {\em Basic} track were black. A smaller percent of black students graduated from DC public school in that year than in 1949. Also, black schools were often overcrowded and offered less curriculum options to students \cite{Findings1965,VoteNo}. In the summer of 1966, Haynes was elected president of the school board, and the board ordered the school administration to dismantle the track system and replace it with new classroom methods \cite{SBPres}.

\begin{figure}
\centering
\includegraphics[height=160mm]{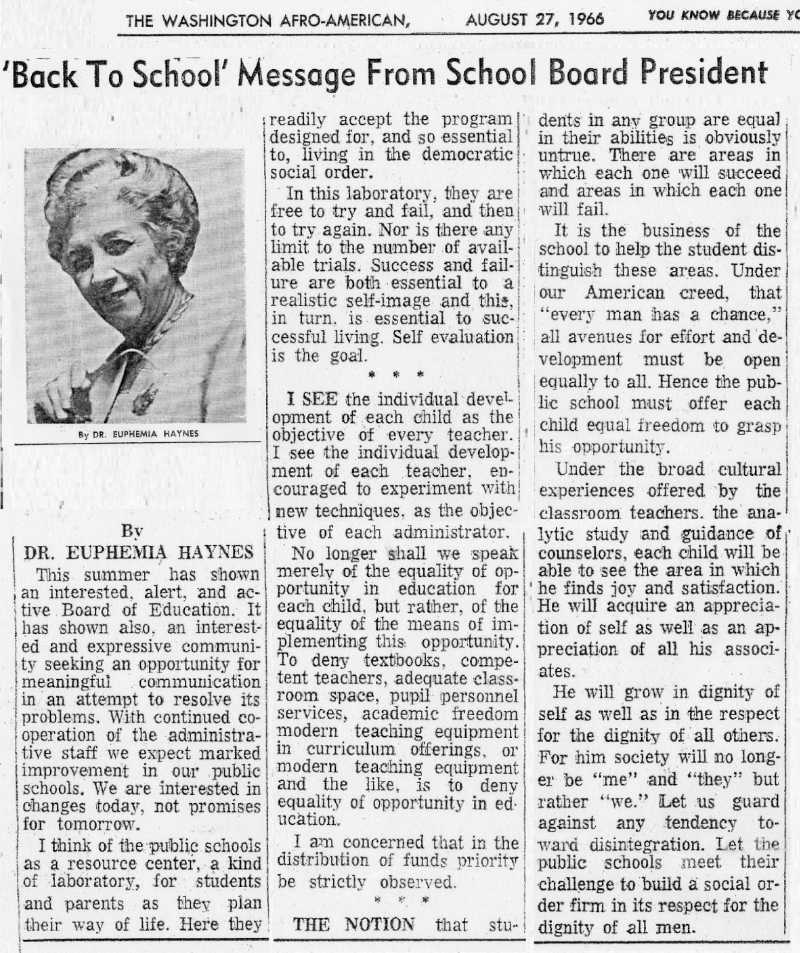}\\
\caption{As the new school board president, Haynes wrote that students ``are free to try and fail, and then try again.''}
\end{figure}

The years of debate culminated in the Hobson v Hansen District Court case that was filed by civil rights activist Julius Hobson in 1966 with the decision handed down on June 21, 1967. Circuit Judge J.~Skelly Wright concluded that the superintendent and school board of DC unconstitutionally deprived black students and poor students their right of an equal education. Among its findings, the court found many examples of discrimination: the school's racial and socially homogeneous make-up was damaging to students; a quota system of four blacks on the nine-member school board kept blacks in the minority even with a public school population made up of 90\% blacks; schools that were predominantly black were overcrowded and underfunded; allowances were given for affluent white students to leave neighborhood schools that were predominantly black; and the track system was based on aptitude tests that were standardized towards middle class white students.   Among the remedies cited, the school district was to abolish the track system and the optimal zones that placed students in neighborhood schools.  Busing was to be offered to students wishing to leave overcrowded schools \cite{Hobson, HvH}.

In a letter from Hansen in December 1966, Superintendent Hansen proposed the option of giving the board the freedom to reappoint him to another term or seek someone new. In March of 1967, Haynes proposed to the board to open a search for a superintendent \cite{Dec66, March67}. Following a majority vote of the board not to appeal the Hobson v Hansen case, on July 3, 1967, Hansen announced his retirement effective July 31 of that year. He wrote that he was unwilling to give school control concerning tracking to Federal sources or the {\em Washington Post}. He had also expressed his belief that ``the predominance of the blacks in Washington population and in the leadership of that population will give them an influence in American life that is out of balance with the eleven per cent proportion to the Nation's population'' \cite{Hansen}.

After the court decision, many community members supported the work of Haynes who many viewed as the track system's most forceful opponent. Both the Passow study and the court decision supported the claims she had forcefully argued since 1963. However, the changes called for also caused in-fighting among school board members, and in an unprecedented move, for the first time, a seated board president was not nominated for a second term. On July 1, 1967, Haynes was deposed as president, and another board member was named in her place.  Haynes remained on the board until 1968 when she chose to retire \cite{HLBiographicalNotes, June27, July1, July18}. It is worth noting that Julius Hobson was elected to the school board in 1968 \cite{VicHob}. For Haynes' service as board president, she was recognized for bringing the DC school system forward further in one year than what had been done in the previous four. Relations were improved between the board and teachers by permitting teachers to select their representative through collective bargaining, and parent and community involvement had been increased with the board's greater transparency on issues. The Passow study and court decision highlighted work that needed to still be done, but Haynes had directed the public schools of DC in a forward direction \cite{Nomination}.

\begin{figure}
\centering
\includegraphics[height=60mm]{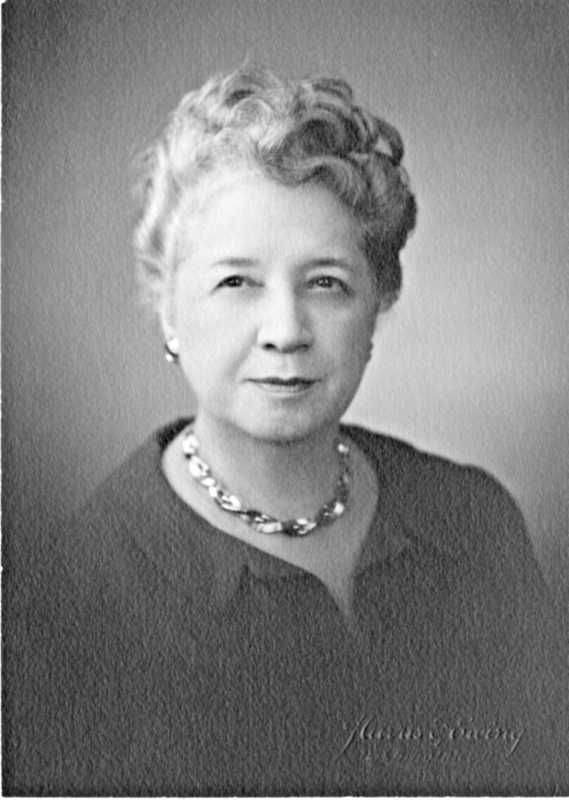}\\
\caption{Haynes was appointed to the DC school board in 1960.  \cite{AppointToSchoolBoard}}
\end{figure}

\section{Thoughts on Gender and Race}

Haynes, with hard work and determination, excelled despite many obstacles for both women and African Americans. She encouraged more women to seek new opportunities: ``The increasing demands by industry for workers together with the change in the nature of household tasks has resulted in a great change in the status of the woman in the home...household tasks ...far simpler and less time consuming than they were a generation ago...If the home of today is to retrieve the blessings and the joys of the traditional home and yet profit by the advantages of the time and labor saving devices of the modern age, it must expand its horizon far beyond the four walls of the family dwelling...her spiritual leadership is not confined to her immediate family, but is felt in her community. Only by influencing the community can she successfully mould her family'' \cite{Home}. After advances made by women during WWII, Haynes also encouraged more women to answer callings to mathematics and science: ``After this war the results of the researches in physics will transform our world.  Whole new industries present a vision of future accomplishment that promises almost miraculous changes in our way of living. This introduction of more physical techniques into industry means a greater demand for physicists.  The industries of tomorrow demand a higher level of mathematical training than those of today. Before the close of this war, women will have established themselves in the technical world.  Tests have shown that women have mechanical ability but have lacked the opportunity to learn about machines... Thousands of women have demonstrated their ability. They have proved that it is worthwhile to get ready for work in science and in industry'' \cite{EHMath}.

Still a greater devotion of Haynes' was to address discrimination based on race. In looking at the world history of slavery, she noted that prior to slavery in America, most masters and slaves were of the same race. The change to different races in America caused an association of inferior status based on the race of a person \cite{EHrace}. Haynes experienced personal discrimination due to race early in life as a student of a segregated school system. Later in life, she also faced discrimination from the shade of her skin.  In a letter of resignation from the Catholic Interracial Council, she explained her reason: During a meeting of the council, they chose representatives to meet with the archbishop during an upcoming reception. They desired the chairman, the secretary and one ``noticeably colored person.'' She stated, `` I am the only colored person on the steering committee of the council, but I do not qualify as `noticeably colored.' A young woman, a former student of mine was selected....I trust that serious consideration will be given by the group to the matter of their own prejudices and that there may soon be a flowering of true Catholic charity which embraces even those Negroes who are not wholly black'' \cite{EHlight}.

Haynes noted that children were not born with prejudices and that parents needed to help by (1) making sure children understand that with current science, people of every country are also our neighbors, (2) helping children to appreciate others with interactions, and (3) getting children to work in groups within their community \cite{Parent}. Citing her faith, she noted ``...race prejudice is not found among children; that among artists, scholars, among those who have achieved something which they have no fear of losing, there seems to be a tendency to lay it aside completely. Because of race prejudice the Negro must endure many forms of social disadvantage...no one of these social ills is half so damaging to him as is the philosophy of life, the type of character he must necessarily develop as the result of racial intolerance and discrimination ... According to the Catholic Doctrine all men are members of one family ... created by God in His own image and likeness... my conviction that a living faith in these fundamental teachings of the Catholic Church and racial friction are incongruous...''\cite{EHCath}.

Finally Hayes tied these thoughts to mathematics and to science: ``...mathematics is an essential factor in cultural integration. ...In whatever corner of the world they may find themselves, mathematicians, like all scientists, are bonded together by a universal desire to understand life.  Cooperation is natural, it is easy, it is necessary in the all-out effort of science to establish truth. As has frequently been pointed out science is international rather than national. The scientific exploration of the universe is the right of all men''\cite{EHmath}.

\section{Conclusion}

Haynes died on July 25, 1980, two years after her husband's passing. The couple had no children, and upon her death, \$700,000 was willed to The Catholic University of America, where a chair in the department of education was named in her honor \cite{HLBiographicalNotes}. It was not until 2001 that Haynes was recognized by the mathematics community as the first African American woman to earn a doctorate in mathematics. Patricia Kenschaft, who conducted early research on black women in mathematics, believed this was due in part to racial segregation both in the education system and in the professional mathematical community during Haynes' life \cite{change}.

Haynes lived life based on her ideals. She saw the need to always take time to fill oneself with knowledge and to then give of oneself in service to others \cite{ValSpeech, SmithLetters, RooseveltHigh}. Haynes spoke of her faith which taught her that everyone ``has a dignity that must be preserved.'' She spoke of the difficulty of achieving goals, but also spoke of the inner strength that her faith gave her: ``...consolation in sorrow...darkness into light...replace disappointment and affliction by peace; fear and hysteria by courage and hope'' \cite{Faith}.

The battles Haynes so gallantly fought are still with us today. The problems exposed in Hobson v Hanson related to neighborhood schools which maintained unequal public education related to prosperous and poorer neighborhoods is still present today.  State funding for public education does not ensure equal assess to quality education for all.  Today, we see additional funding needs addressed with local school referendums, and thus wealthier regions can afford to vote to give more to schools while other regions may not have such resources \cite{WPR}. In addition, tracking, ability grouping, and gifted and talented programs still face similar issues.  While schools look for ways to best teach diverse groups of students, such systems of separation still often occur along racial lines and income levels. While advanced classes are found to academically benefit the students who take them, students in lower-level classes are found to suffer. Current suggestions for improvement include requiring all students to take tests to qualify for advanced tracks, thus eliminating the common practice of such tests being taken predominantly by students whose parents push for such options.  Another suggestion, similar to Haynes' ideas, is to remove all barriers to enter advanced classes and let students self-select instead of basing entry on exams \cite{Atlantic}.

Haynes devoted much of her life to her own education and to the education of public school students of Washington DC.  In examining the battles she fought to improve educational opportunities and in reflecting on the progress that still must be made, we are wise to consider her words in an address delivered to educators in 1945: ``The concept of limit is merely an expression mathematical form of an ever receding goal of perfection for which man yearns and for which he strives, yet never attains. With each new approximation, he is merely closer'' \cite{HaynesLogicQuote}. Haynes knew that  some goals are difficult to achieve and sometimes are not reached.  That did not mean they should not be attempted. As we continue to face the challenges of equal opportunities for all, we are wise to remember that we may not reach the goal, but we can each work to move the world closer to the ``goal of perfection.''


\def\b{\bibitem}

 { \footnotesize

\newpage

\noindent  The images in this paper are courtesy of the following sources:
 \newline Figure 1 Photo of Euphemia Haynes recieving Lady of Year Award from the Oldest Inhabitants, The Evening Star, April 24, 1967, Box 3, Folder 10, Haynes-Lofton Family Papers, ACUA, The Catholic University of America, Washington DC.
 \newline Figure 2 Courtesy of 1860 United States Census of Washington DC., database with images, FamilySearch, \url{https://familysearch.org/ark:/61903/1:1:MCV9-9G3} (accessed 10 May 2016) Mary A Day (Schedual I - Free Inhabitants; pg. 42; dwelling 258; family 276; enum. June 19, 1860).
 \newline Figure 3 Photo of Euphemia Lofton Haynes, Box 67, Folder 24, Haynes-Lofton Family Papers, ACUA, The Catholic University of America, Washington DC.
 \newline Figure 4 Photo of Euphemia and Harold Haynes, Box 3, Folder 10, Haynes-Lofton Family Papers, ACUA, The Catholic University of America, Washington DC.
 \newline Figure 5 Photo of Aubrey E. Landry, Box 7, File 3, Catholic University of America Photograph Collection ca. 1886-2000, ACUA, The Catholic University of America, Washington DC.
 \newline Figure 6 Euphemia Lofton Haynes Doctoral Thesis cover page, ACUA, The Catholic University of America, Washington DC.
 \newline Figure 7 Photo of Euphemia Haynes Graduation from The Catholic University of America, Box 67, Folder 21, Haynes-Lofton Family Papers, ACUA, The Catholic University of America, Washington DC.
 \newline Figure 8 Photo of Euphemia Haynes and unidentified man, Box 67, Folder 17, Haynes-Lofton Family Papers, ACUA, The Catholic University of America, Washington DC.
 \newline Figure 9 USO Photo of Euphemia Haynes, Box 67, Folder 17, Haynes-Lofton Family Papers, ACUA, The Catholic University of America, Washington DC.
 \newline Figure 10 'Back to School' Message from School Board President, The Washington Afro-American, Augusts 27, 1966, Box 50, Folder 3, Haynes-Lofton Family Papers, ACUA, The Catholic University of America, Washington DC.
 \newline Figure 11 Photo of Euphemia Lofton Haynes, Box 67, Folder 23, Haynes-Lofton Family Papers, ACUA, The Catholic University of America, Washington DC.

 \newpage

\begin{figure}
\centering
\includegraphics[height=60mm]{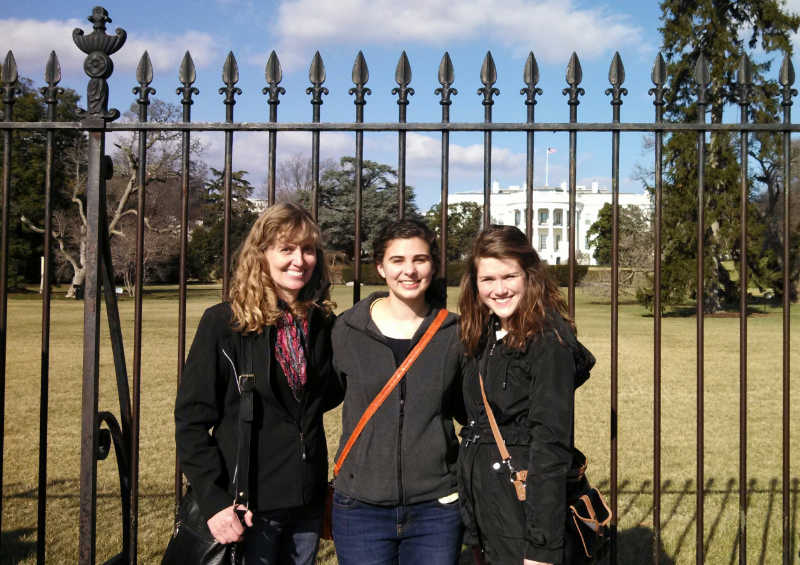}\\
\caption{Kelly, Zoroufy and Shinner in Washington DC during trip to the Archive Library of The Catholic University of America.}
\end{figure}

 Funding for travel to Washington DC to study the archive materials of the Haynes-Lofton Family, located at The Catholic University of America Archive Library, and to obtain copies of needed items was provided by a University of Wisconsin - La Crosse Undergraduate Research Grant.  This grant made it possible to obtain most of the materials used in this paper and also provided us a memorable research experience.

 Appreciation is extended to Paul Kelly, former Archive Librarian at the Catholic University of America, who helped with our work prior to our trip as well and during and after our visit. We also wish to thank Richard M Lofton, step-son of Euphemia Lofton's brother, Joseph Lofton, for talking with us about this project in 2014.  Appreciation is also given to Richard Breaux, Tracy Helixon, Greg Wegner, and Barbara Weiss for support, suggestions and feedback on this paper.

The authors dedicate this paper to public school teachers across the country who continue to work towards Euphemia's goals of quality and equal educational opportunities for all students. Although this goal may never be fully realized, each step can get us closer to that goal.


\end{document}